\documentclass{amsart}
\input epsf

\usepackage{amsfonts,amsthm,amsmath,amssymb,latexsym}
\usepackage{graphicx,color}
\usepackage[all]{xy}

\begin{document}

\author{Dragomir \v Sari\' c}
\thanks{This research was partially supported by National Science Foundation grant DMS 1102440.}

\address{Department of Mathematics, Queens College of CUNY,
65-30 Kissena Blvd., Flushing, NY 11367}
\email{Dragomir.Saric@qc.cuny.edu}

\address{Mathematics PhD. Program, The CUNY Graduate Center, 365 Fifth Avenue, New York, NY 10016-4309}

\theoremstyle{definition}

 \newtheorem{definition}{Definition}[section]
 \newtheorem{remark}[definition]{Remark}
 \newtheorem{example}[definition]{Example}

\newtheorem*{notation}{Notation}

\theoremstyle{plain}

 \newtheorem{proposition}[definition]{Proposition}
 \newtheorem{theorem}[definition]{Theorem}
 \newtheorem{corollary}[definition]{Corollary}
 \newtheorem{lemma}[definition]{Lemma}

\newcommand{\eps}{\varepsilon}
\newcommand{\G}{\Gamma}
\newcommand{\g}{\gamma}
\newcommand{\D}{\Delta}
\renewcommand{\d}{\delta}
\newcommand{\dist}{\mathrm{dist}}
\newcommand{\m}{\mathrm{mod}}

\title[Asymptotically conformal deformations]{Fenchel-Nielsen coordinates for asymptotically conformal deformations}


\date{\today}

\maketitle

\begin{abstract}
Let $X$ be an infinite hyperbolic surface endowed with an upper bounded geodesic pants decomposition.  Alessandrini, Liu, Papadopoulos, Su and Sun \cite{ALPSS}, \cite{ALPS} parametrized the quasiconformal Teichm\"uller space $T_{qc}(X)$ and the length spectrum Teichm\"uller space $T_{ls}(X)$ using the Fenchel-Nielsen coordinates. A quasiconformal map $f:X\to Y$ is said to be {\it asymptotically conformal} if its Beltrami coefficient $\mu =\bar{\partial}f/\partial f$ converges to zero at infinity. The space of all asymptotically conformal maps up to homotopy and post-composition by conformal maps is called ``little'' Teichm\"uller space $T_0(X)$. We find a parametrization of $T_0(X)$ using the Fenchel-Nielsen coordinates and a parametrization of the closure $\overline{T_0(X)}$ of $T_0(X)$ in the length spectrum metric. We also prove that the quotients $AT(X)=T_{qc}(X)/T_0(X)$,  $T_{ls}(X)/\overline{T_{qc}(X)}$ and $T_{ls}(X)/\overline{T_0(X)}$ are contractible in the Teichm\"uller metric and the length spectrum metric, respectively. Finally, we show that the Wolpert's lemma on the lengths of simple closed geodesics under quasiconformal maps is not sharp.
\end{abstract}

\section{Introduction}

Let $X$ be a complete hyperbolic surface with infinitely generated fundamental group whose boundary, if non-empty, consists of simple closed geodesics. Let $\mathcal{P}=\{\alpha_n\}$ be a geodesic pants decomposition of $X$ by simple closed geodesics $\alpha_n$, where some $\alpha_n$ can have zero length i.e. $\alpha_n$ can be a puncture. In addition, assume that there exists a fixed constant $P>0$ such that
$$
l_{\alpha_n}(X)\leq P
$$
for all $n$, where $l_{\alpha_n}(X)$ is the length of the geodesic $\alpha_n$ for the hyperbolic metric of $X$.
We say that the pants decomposition $\mathcal{P}=\{\alpha_n\}$ is {\it upper bounded}. 

The {\it quasiconformal} Teichm\"uller space $T_{qc}(X)$ consists of all quasiconformal maps $f:X\to Y$ up to post-composition by isometries and up to  homotopies. 
The {\it length spectrum} Teichm\"uller space $T_{ls}(X)$ 
consists of all homeomorphisms $h:X\to Y$ up to post-composition by isometries and up to homotopies such that
$$
L(X,Y):=\sup_{\beta}\max\{ \frac{l_{\beta}(X )}{l_{\beta}(Y )},\frac{l_{\beta}(Y )}{l_{\beta} (X)}\}<\infty ,
$$
where the supremum is over all simple closed curves $\beta$ on $X$,
and where $l_{\beta}(X ),l_{\beta}(Y )$ are the lengths of the geodesic representatives of $\beta$ on $X,Y$, respectively.  Wolpert (cf. \cite{Wol1}) proved that if $f:X\to Y$ is a K-quasiconformal map between hyperbolic surfaces then for each non-trivial simple closed curve $\beta$ on $X$ we have
\begin{equation}
\label{eq:Wolpert}
\frac{1}{K}l_{\beta}(X)\leq l_{\beta}(Y)\leq Kl_{\beta}(X)
\end{equation}
which gives $T_{qc}(X)\subset T_{ls}(X)$.

Let $\{ (l_n(Y),t_n(Y))\}_n$ denote the Fenchel-Nielsen coordinates for a marked hyperbolic surface $Y$ corresponding to the pants decomposition $\mathcal{P}=\{\alpha_n\}$, where $(l_n(Y),t_n(Y))$ is defined for all $n$ with $l_{\alpha_n}(X)>0$ and $\alpha_n$ not a boundary component. If $\alpha_n$ is a boundary geodesic, then the twist $t_n(Y)$ is not defined. Alessandrini, Liu, Papadopoulos, Su and Sun \cite{ALPSS} proved that a sequence $\{ (l_n,t_n)\}_n$ with $l_n>0$ and $t_n\in\mathbb{R}$ defined for all $n$ with $l_{\alpha_n}(X)>0$ are the Fenchel-Nielsen coordinates of $Y\in T_{qc}(X)$ if and only if there exists $M>0$ such that
$$
|\log l_n/l_n(X)|\leq M
$$
and
$$
|t_n-t_n(X)|\leq M.
$$
Moreover, the Fenchel-Nielsen map given by $Y\mapsto \{ (\log l_n(Y), t_n(Y))\}_n$ is a locally bi-Lipschitz parametrization of $T_{qc}(X)$ equipped with the Teichm\"uller metric by the space of bounded sequences $l^{\infty}$ equipped with the supremum norm (cf. \cite{ALPSS}).

A quasiconformal map $f:X\to Y$ is said to be {\it asymptotically conformal} if for every $\epsilon >0$ there exists a compact set $K\subset X$ such that $\|\mu|_{X-K}\|_{\infty}<\epsilon$, where $\mu=\bar{\partial}f/\partial f$ is the Beltrami coefficient of $f$. The ``little'' Teichm\"uller space $T_0(X)$ consists of all asymptotically conformal quasiconformal maps $f:X\to Y$ up to homotopy and post-composition by isometries. Note that $T_0(X)$ is a closed nowhere dense subset of $T_{qc}(X)$. We prove

\vskip .2 cm

\noindent
{\bf Theorem 1.} {\it Let $X$ be an infinite hyperbolic surface equipped with an upper bounded geodesic pants decomposition $\mathcal{P}=\{\alpha_n\}_n$. Then a sequence $\{ (l_n,t_n)\}_n$ with $l_n >0$ and $t_n\in\mathbb{R}$ represents a point in $T_0(X)$ if and only if
\begin{equation}
\label{eq:little-lengths}
|\log\frac{l_n}{l_n(X)}|\to 0
\end{equation}
and 
\begin{equation}
\label{eq:little-twists}
|t_n-t_n(X)|\to 0
\end{equation}
as $n\to\infty$.
}

\vskip .2 cm

\noindent {\bf Remark 1.} The equation (\ref{eq:little-lengths}) is the asymptotically conformal version of the Wolpert's inequality (\ref{eq:Wolpert}).

\vskip .2 cm

The {\it length spectrum distance} $d_{ls}$ between two marked surfaces $X$ and $Y$ is defined by
$$
d_{ls}(X,Y)=\frac{1}{2}\log L(X,Y).
$$
Shiga \cite{Shiga} started the study of the length spectrum distance $d_{ls}$ on $T_{qc}(X)$ for infinite surfaces and he established that $d_{ls}$ is not complete on $T_{qc}(X)$ in general. Alessandrini, Liu, Papadopoulos and Su \cite{ALPS1} proved that $T_{ls}(X)$ is complete in the length spectrum distance $d_{ls}$. 
Moreover, the closure of $T_{qc}(X)$ for the length spectrum distance is nowhere dense in $T_{ls}(X)$ (cf. \cite{ALPS}). Alessandrini, Liu, Papadopoulos and Su \cite{ALPS} proved that 
a sequence $\{ (l_n,t_n)\}_n$ with $l_n>0$ and $t_n\in\mathbb{R}$ represents a point in $T_{ls}(X)$ if and only if there exists $M>0$ such that 
\begin{equation}
\label{eq:ls-lengths}
|\log\frac{l_n}{l_n(X)}|\leq M
\end{equation}
and 
\begin{equation}
\label{eq:ls-twists}
\frac{|t_n-t_n(X)|}{\max\{ 1,|\log l_n(X)|\} }\leq M
\end{equation}
for all $n$.

The {\it normalized Fenchel-Nielsen map} defined by
\begin{equation*}
F(Y)=\Big{\{} (\log\frac{l_{n}(Y)}{l_{n}(X)},\frac{t_{n}(Y)-t_{n}(X)}{\max\{ 1,|\log l_{n}(X)|\}})\Big{\}}_n
\end{equation*}
is a locally bi-Lipschitz of $T_{ls}(X)$ equipped with the length-spectrum metric onto $l^{\infty}$ equipped with the supremum norm (cf. \cite{Sar}). In addition, a marked surface $Y$ is in the closure $\overline{T_{qc}(X)}$ of $T_{qc}(X)$ in the length spectrum metric if and only if there exists $M>0$ such that the Fenchel-Nielsen coordinates $\{ (l_n(Y),t_n(Y))\}$ of $Y$ satisfy (cf. \cite{Sar})
$$
|\log\frac{l_n(Y)}{l_n(X)}|\leq M
$$
and
$$
| t_n(Y)-t_n(X)|=o(|\log l_n(X)|)
$$
as $l_n(X)\to 0$. (If a sequence $a_n\to\infty$ then $o(a_n)$ satisfies $o(a_n)/a_n\to 0$.) We prove

\vskip .2 cm

\noindent {\bf Theorem 2.} {\it 
Let $X$ be an infinite complete hyperbolic surface with geodesic boundary, if any,  that has an upper bounded geodesic pants decomposition. Let  $\overline{T_0(X)}$ be the closure in the length spectrum metric of the little Teichm\"uller space $T_0(X)$. Then $Y\in \overline{T_0(X)}$ if and only if its Fenchel-Nielsen coordinates satisfy
\begin{equation}
\label{eq:ls-l0}
|\log\frac{l_n(Y)}{l_n(X)}|\to 0
\end{equation}
and
\begin{equation}
\label{eq:ls-t0}
\frac{|t_n(Y)-t_n(X)|}{\max\{ 1,|\log l_n(X)|\}}\to 0
\end{equation}
as $n\to\infty$.
}

\vskip .2 cm

We obtain the following corollary

\vskip .2 cm

\noindent {\bf Corollary 1.} {\it  
Let $X$ be an infinite complete hyperbolic surface with geodesic boundary, if any, that has an upper bounded geodesic pants decomposition. Then the asymptotic Teichm\"uller space $AT(X)=T_{qc}(X)/T_0(X)$ is contractible in the Teichm\"uller metric. Moreover, quotients $T_{ls}(X)/\overline{T_{qc}(X)}$ and $T_{ls}(X)/\overline{T_0(X)}$ are contractible in the length spectrum metric.
}

\vskip .2 cm

Fletcher \cite{Fletcher1}, \cite{Fletcher} proved that $AT(X)$ is locally bi-Lipschitz to the space $l^{\infty}/c_0$  for the {\it unreduced} Teichm\"uller spaces, where $l^{\infty}$ is the space of bounded sequences and $c_0$ is the subspace of sequences converging to zero. The above parameterizations and \cite{Sar} give that $AT(X)$,  $T_{ls}(X)/\overline{T_{qc}(X)}$ and $T_{ls}(X)/\overline{T_0(X)}$ are globally parametrized and locally bi-Lipschitz to $l^{\infty}/c_0$ where the spaces are defined in the {\it reduced} sense, i.e. homotopies are allowed to move points on the boundary geodesics. 

\vskip .2 cm

For each $\alpha_n$ in the pants decomposition $\mathcal{P}$ of $X$ with $l_{\alpha_n}(X)>0$, let $P_n^1$ and $P_n^2$ be the two pairs of pants of $\mathcal{P}$ that have $\alpha_n$ on their boundary with possibly $P_n^1=P_n^2$.
Define $\gamma_n$ to be a shortest closed geodesic in $P_n^1\cup P_n^2$ intersecting $\alpha_n$ in either two or one point. 
We have

\vskip .2 cm

\noindent {\bf Proposition 1.} {\it Let $X$ be an infinite complete hyperbolic surface whose boundary, if any, consists of simple closed geodesics. Let $f:X\to Y$ be a quasiconformal map. Then
$$
\sup_n |l_{\gamma_n}(Y)-l_{\gamma_n}(X)|<\infty .
$$
}

\vskip .2 cm 

\noindent
{\bf Remark 2.} The above proposition proves that the Wolpert's inequality (\ref{eq:Wolpert}) is not sharp if there exists a subsequence $\gamma_{n_k}$ with $l_{\gamma_{n_k}}(X)\to\infty$ as $n_k\to\infty$. This is the case when $l_{\alpha_{n_k}}(X)\to 0$ as $n_k\to\infty$.

\section{The Fenchel-Nielsen coordinates for $T_0(X)$}

Let $X$ be an infinite surface (i.e. the fundamental group of $X$ is infinitely generated) that is equipped with a complete hyperbolic metric such that the boundary of $X$, if non-empty, consists of simple closed geodesics. Assume $X$ has a generalized geodesic pants decomposition $\mathcal{P}=\{\alpha_n\}_n$ ({\it generalized pants decomposition} means that some of the boundary components $\alpha_n$ of the pants could be punctures) such that there exists a constant $P>0$ 
with
$$
l_{\alpha_n}(X)\leq P.
$$
Here $l_{\alpha_n}(X)$ stand for the length of the geodesic $\alpha_n$ in the hyperbolic metric of $X$. If $\alpha_n$ is a puncture, then $l_{\alpha_n}(X)=0$. We say that the geodesic pants decomposition $\mathcal{P} =\{\alpha_n\}$ is {\it upper bounded.}

A quasiconformal map $f:X\to Y$ is {\it asymptotically conformal} if for every $\epsilon >0$ there exists a compact set $K_{\epsilon}\subset X$ such that the Beltrami coefficient $\mu =\bar{\partial}f/\partial f$ satisfies $\|\mu |_{X-K_{\epsilon}}\|_{\infty}<\epsilon$. Note that the compact set $K_{\epsilon}$ might contain finitely many  boundary geodesics of $X$. The quasiconformal Teichm\"uller space $T_{qc}(X)$ consists of all quasiconformal maps $f:X\to Y$ up to post-composition by hyperbolic isometries of the image surface and up to homotopy. Note that the homotopies are fixing boundary geodesics set-wise but they can move points on the boundary. The quasiconformal Teichm\"uller space $T_{qc}(X)$ is also known as the {\it reduced} Teichm\"uller space. The ``little'' Teichm\"uller space $T_0(X)$ consists of all asymptotically conformal quasiconformal maps $f:X\to Y$ up to post-composition by hyperbolic isometries of the image surface and up to homotopy. The homotopies set-wise fix boundary geodesics of $X$.

Given a homeomorphism $h:X\to Y$, a closed geodesic $\beta$ on $X$ is mapped onto a closed curve $h(\beta )$ on $Y$. We denote by $l_{\beta}(Y)$ the length of the unique closed geodesic in $Y$ homotopic to $\beta$. If $\alpha_n$ is a geodesic of the pants decomposition $\mathcal{P}$ of $X$, denote $l_n(Y)=l_{\alpha_n}(Y)$ the length part of the Fenchel-Nielsen coordinates for $Y$. If $\alpha_n$ is a boundary geodesic then the twist $t_n$ is not defined. If $\alpha_n$ is interior geodesic of $X$, then we define the twist parameter $t_n(Y)$ as the twist parameter of a compact subsurface of $X$ with geodesic boundary which contains $\alpha_n$ in its interior. It is irrelevant which subsurface we choose (cf. \cite{ALPS}). Therefore we obtained the Fenchel-Nielsen coordinates $\{ (l_n(Y),t_n(Y)\}$ for any marked surface $Y$, where $(l_n,t_n)$ is not specified when $l_{\alpha_n}(X)=0$ and $t_n$ is not specified when $\alpha_n$ is a boundary geodesic.

\vskip .2 cm

\noindent {\it Proof of Theorem 1.} 
Assume first that $\{ (l_n,t_n)\}$ is given that satisfies (\ref{eq:ls-lengths}) and (\ref{eq:ls-twists}). 
Denote by $Y$ a hyperbolic surface homeomorphic to $X$ with the Fenchel-Nielsen coordinates $\{ (l_n,t_n)\}$ (cf. \cite{ALPS}). Note that $Y$ is constructed by gluing geodesic pairs of pants according to the twists $t_n$. The shape of the geodesic pants is determined by the lengths $l_n$. Following \cite{ALPS}, we construct a quasiconformal map $f:X\to Y'$ such that the surface $Y'$ has Fenchel-Nielsen coordinates $\{ l_n,\frac{l_n}{l_n(X)}t_n(X)\}$.
It is immediate from \cite{ALPS} and Bishop \cite{Bishop} that the quasiconformal map $f:X\to Y'$ is aymptotically conformal by the condition (\ref{eq:ls-lengths}). It remains to prove that there is an asymptotically conformal map $g:Y'\to Y$. Then $g\circ f:X\to Y$ is an asymptotically conformal map. Note that each geodesic boundary $\alpha_n$ in the interior of $Y'$ has a collar neighborhood of width bounded below by a positive constant. Thus the twisting by the amount $t_n-\frac{l_n}{l_n(X)}t_n(X)$ along the family $\alpha_n$ is homotopic to a quasiconformal map obtained by an explicit construction in each collar neighborhood (cf. Wolpert \cite{Wol1}) because $|t_n-\frac{l_n}{l_n(X)}t_n(X)|$ is bounded. Moreover, conditions (\ref{eq:ls-lengths}) and (\ref{eq:ls-twists}) imply that  $|t_n-\frac{l_n}{l_n(X)}t_n(X)|\to 0$ as $n\to\infty$ and the construction in \cite{Wol1} implies that the map $g$ can be chosen to be asymptotically conformal.

It remains to prove that if $f:X\to Y$ is asymptotically conformal $K$-quasinformal map, then its Fenchel-Nielsen coordinates $\{ (l_n(Y),t_n(Y))\}_n$ satisfy (\ref{eq:ls-lengths}) and (\ref{eq:ls-twists}). 
Let $\mathbb{H}$ be the upper half-plane model of the hyperbolic plane.
We consider a sequence of universal coverings $\pi_n^1:\mathbb{H}\to X$ and $\pi_n^2:\mathbb{H}\to Y$ such that positive half of the $y$ axis covers the geodesic $\alpha_n\subset X$ and the geodesic homotopic to $f(\alpha_n)\subset Y$, respectively. Let $\tilde{f}_n:\mathbb{H}\to\mathbb{H}$ be the lift of $f:X\to Y$ that fixes $0$, $1$ and $\infty$. Denote by $A_n^1$ and $A_n^2$ the primitive hyperbolic elements of the covering groups for $\pi_n^1$ and $\pi_n^2$ for the curves $\alpha_n\subset X$ and $f(\alpha_n)\subset Y$, respectively. Then
$$
\tilde{f}_n\circ A_n^1\circ \tilde{f}_n^{-1}=A_n^2
$$
and the translation lengths of $A_n^1$ and $A_n^2$ are $l_n(X)$ and $l_n(Y)$, respectively.
Then $K\cdot P$ is an upper bound for the lengths of the pants decomposition $\{ f(\alpha_n)\}$ on $Y$.

Let $a_n\in\mathbb{N}$ be such that $P\leq a_nl_n(X)\leq 2P$. To prove (\ref{eq:ls-lengths}) assume that it is false and find a contradiction. Namely there exists a subsequence $l_{n_k}(Y)$ and $C>0$ such that 
\begin{equation}
\label{eq:ls-conv}
|\log\frac{l_{n_k}(Y)}{l_{n_k}(X)}|\geq C>0.
\end{equation}
The sequence of quasiconformal maps $\tilde{f}_{n_k}:\mathbb{H}\to\mathbb{H}$ has the same quasiconformal constant and it fixes $0$, $1$ and $\infty$. Therefore it has a subsequence that converges uniformly on compact subsets of $\mathbb{H}$ to a quasiconformal map $\tilde{f}:\mathbb{H}\to\mathbb{H}$ and for simplicity of notation denote this subsequence by $\tilde{f}_{n_k}$. The sequence of Beltrami coefficients $\mu_{n_k}=\bar{\partial}\tilde{f}_{n_k}/\partial \tilde{f}_{n_k}$ converges to zero pointwise a.e. on $\mathbb{H}$ as $n_k\to\infty$ because $\tilde{f}_{n_k}$ is a lift of an asymptotically conformal quasiconformal map $f$ such that $\alpha_{n_k}$ lifts to the positive $y$ axis. Thus any fixed compact subset of $\mathbb{H}$ covers parts of the surface $X$ which leave every compact subset of $X$ which forces convergence to zero of $\mu_{n_k}$. Therefore $\tilde{f}$ is a conformal map which fixes $0$, $1$ and $\infty$, and therefore $\tilde{f}=id$.

Consider the hyperbolic isometries $(A_{n_k}^1)^{a_{n_k}}$ and $(A_{n_k}^2)^{a_{n_k}}$ which fix $0$ and $\infty$. They have translation lengths between $\frac{1}{K}P$ and $2KP$. Therefore there exist subsequences of $(A_{n_k}^1)^{a_{n_k}}$ and $(A_{n_k}^2)^{a_{n_k}}$ which pointwise converge to hyperbolic isometries $A^1$ and $A^2$ both different from the identity. Since
$$
\tilde{f}_{n_k}\circ (A_{n_k}^1)^{a_{n_k}}\circ \tilde{f}_{n_k}^{-1}=(A_{n_k}^2)^{a_{n_k}}
$$
and $\tilde{f}_{n_k}\to id$ as $n_k\to\infty$, we have
$$A^1=A^2.$$
In particular they have the same translation lengths. This implies that the ratio of the translation length of $(A_{n_k}^1)^{a_{n_k}}$ to the translation length of  
$(A_{n_k}^2)^{a_{n_k}}$ converges to $1$, namely
$$
\frac{a_{n_k}l_{n_k}(Y)}{a_{n_k}l_{n_k}(X)}=\frac{l_{n_k}(Y)}{l_{n_k}(X)}\to 1
$$
as $n_k\to\infty$. This contradicts (\ref{eq:ls-conv}) and establishes (\ref{eq:ls-lengths}).

Let $Y'$ be the marked hyperbolic surface whose Fenchel-Nielsen coordinates are $\{ (l_n(Y),t_n(X))\}$. There exists a marking map $g:X\to Y'$ that is asymptotically conformal quasiconformal map by (\ref{eq:ls-lengths}) and (\ref{eq:ls-twists}), and by the first part of this proof. Then the map $f_1=f\circ g^{-1}:Y'\to Y$ is an asymptotically conformal $K_1$-quasiconformal map. To prove (\ref{eq:ls-twists}), assume on the contrary that there exists a subsequence $t_{n_k}(Y)$ and $C>0$ such that
\begin{equation}
\label{eq:ls-sub>0}
|t_{n_k}(Y)-t_{n_k}(X)|=|t_{n_k}(Y)-t_{n_k}(Y')|\geq C>0
\end{equation}
for all $n_k$. 

We seek a contradiction with (\ref{eq:ls-sub>0}).
For $\alpha_{n}$ and $g:X\to Y'$, the geodesic on $Y'$ homotopic to $g(\alpha_{n})$ is denoted by $g(\alpha_{n})$ for simplicity. Denote by $\mathcal{P}$ the corresponding geodesic pants decomposition of $Y'$. Let $\gamma_{n_k}$ be a shortest simple closed geodesic in $Y'$ which intersects $g(\alpha_{n_k})$ and is contained in the union of two geodesic pairs of pants of $\mathcal{P}$ with the geodesic $g(\alpha_{n_k})$ on their boundaries. We fix universal coverings $\pi_{n_k}^1:\mathbb{H}\to Y'$ and $\pi_{n_k}^2:\mathbb{H}\to Y$ such that lifts of the geodesic  $g(\alpha_{n_k})$ of $Y'$ and geodesic of $Y$ homotopic to $f(\alpha_{n_k})$ contain the positive $y$ axis, and that a lift $\tilde{\gamma}_{n_k}$ of $\gamma_{n_k}$ has one endpoint at $1\in\mathbb{R}$.
The map $f_1:Y'\to Y$ represents twisting along the family $\mathcal{P}$. Let $\tilde{\mathcal{P}}$ be the lift to $\mathbb{H}$ of the pants decomposition $\mathcal{P}$ of $Y'$.
 We choose a lift $\tilde{f}_{n_k}:\mathbb{H}\to\mathbb{H}$ of $f_1$ that fixes $0$ and $\infty$, and that is further normalized such that $\tilde{f}_{n_k}$ fixes the ideal points on $\bar{\mathbb{R}}=\partial\mathbb{H}$ of the component of $\mathbb{H}-\tilde{\mathcal{P}}$ on the left side of the positive $y$ axis and adjacent to the positive $y$ axis.

Note that $\tilde{\gamma}_{n_k}$ intersect only the lifts of $g(\alpha_{n_k})$ from all the geodesics in the pants decomposition $\mathcal{P}$ of $Y'$. Therefore the twisting along $\tilde{\gamma}_{n_k}$ is in the same direction (cf. Alessandrini, Liu, Papadopoulos and Su \cite{ALPS}). Let $x_{n_k}<0$ be the  other endpoint of $\tilde{\gamma}_{n_k}$. Then there exists $J>1$ such that
$$
\frac{1}{J}<\log cr(x_{n_k},0,1,\infty )=\log \frac{1-x_{n_k}}{-x_{n_k}}=\log (\frac{1}{-x_{n_k}}+1)<J
$$
for all $n$ by the choice of $\alpha_{n_k}$ and $\gamma_{n_k}$, and by the upper bound on the lengths of $\alpha_n$ (cf. \cite{ALPSS}). In particular, $x_{n_k}$ is contained in a compact subset of $\mathbb{R}$. Assume without loss of generality that $t_n(Y)-t_n(X)>0$.
By the choice of the lift $\tilde{f}_{n_k}$, we have that $\tilde{f}_{n_k}(0)=0$, $\tilde{f}_{n_k}(\infty )=\infty$, $\tilde{f}_{n_k}(1)\geq e^{t_{n_k}(Y)-t_{n_k}(X)}\geq e^C>1$ and $\tilde{f}_{n_k}(x_{n_k})\geq x_{n_k}$.
It follows that
$$
\log cr(\tilde{f}_{n_k}(x_{n_k}),\tilde{f}_{n_k}(0),\tilde{f}_{n_k}(1),\tilde{f}_{n_k}(\infty ))
\geq \log (\frac{e^C}{-x_{n_k}}+1).
$$
Thus there exists $C'>0$ such that
\begin{equation}
\label{eq:cr>C}
\log \frac{cr(\tilde{f}_{n_k}(x_{n_k}),\tilde{f}_{n_k}(0),\tilde{f}_{n_k}(1),\tilde{f}_{n_k}(\infty ))}{cr(x_{n_k},0,1,\infty )}\geq C'>0
\end{equation}
for all $n_k$. 

On the other hand, $\tilde{f}_{n_k}$ fixes $0$ and $\infty$, and $\tilde{f}_{n_k}(1)$ is bounded away from $0$ and $\infty$ since the total twisting along $\tilde{\gamma}_{n_k}$ is Thurston bounded (cf. \cite{Sar2}). We reflect $\tilde{f}_{n_k}$ in $\mathbb{R}$ to obtain a quasiconformal mapping of $\bar{\mathbb{C}}$. Since the Beltrami coefficient of $\tilde{f}_{n_k}$ converges pointwise a.e. to zero, it follows that $\tilde{f}_{n_k}$ converges to the identity uniformly on compact subsets of $\mathbb{C}$. This contradicts (\ref{eq:cr>C}) and proves (\ref{eq:ls-twists}). {\it End of proof of Theorem 1.}

\section{The closure of $T_{0}(X)$}

We consider the closure $\overline{T_0(X)}$ of $T_0(X)$ as a subset of $T_{ls}(X)$ equipped with the length spectrum metric and find its characterization in terms of Fenchel-Nielsen coordinates.

\vskip .2 cm 

\noindent {\it Proof of Theorem 2.}  
Consider a sequence of marked surfaces $Y_k\in T_0(X)$ that converges to $Y\in \overline{T_{0}(X)}$. Let $\{ (l_n(Y_k),t_n(Y_k))\}$ be the Fenchel-Nielsen coordinates of $Y_k$. 
By Theorem 1, we have
$$
|\log \frac{l_n(Y_k)}{l_n(X)}|\to 0
$$
and
$$
|t_n(Y_k)-t_n(X)|\to 0
$$
as $n\to\infty$ for each fixed $k$. Since the normalized Fenchel-Nielsen map is homeomorphism for the length spectrum distance, the above conditions imply that the limit $Y$ satisfies (\ref{eq:ls-l0}) and (\ref{eq:ls-t0}).

Assume now that $Y\in\overline{T_{qc}(X)}$ has Fenchel-Nielsen coordinates satisfying (\ref{eq:ls-l0}) and (\ref{eq:ls-t0}). We need to find a sequence $Y_k\in T_0(X)$ that converges to $Y$ as $k\to\infty$. Let $Y_k$ be defined by the Fenchel-Nielsen coordinates satisfying
$$
l_n(Y_k)=l_n(Y),\ t_n(Y_k)=t_n(Y)
$$
for $n\leq k$, and
$$
l_n(Y_k)=l_n(X),\ t_n(Y_k)=t_n(X)
$$
for $n>k$. Then $Y_k\in T_0(X)$ by Theorem 1. Since 
$l_n(Y_k)=l_n(Y)$ for $n\leq k$ and $\lim_{n\to\infty}\frac{l_n(Y)}{l_n(X)}=1$, we have that $\sup_n |\log \frac{l_n(Y_k)}{l_n(Y}|\to 0$ as $k\to\infty$. Since $t_n(Y_k)=t_n(X)$ for $n\leq k$, since $|t_n(Y_k)-t_n(Y)|=|t_n(X)-t_n(Y)|$ for $n>k$, and since $|t_n(X)-t_n(Y)|=o(|\log l_n(X)|)$ (by $t_n(Y)\in\overline{T_{qc}(X)}$), we get that $Y_k\to Y$ as $k\to\infty$ in the length spectrum metric. 
{\it End of proof of Theorem 2.}

\vskip .2 cm

\noindent {\it Proof of Corollary 1.} The Fenchel-Nielsen coordinates map $T_{qc}(X)$ homeomorphically onto $l^{\infty}$ and they map $T_0(X)$ homeomorphically onto the space $c_0$ of sequences 
that vanish at infinity. Then $AT(X)$ is homeomorphic to $l^{\infty}/c_0$ and thus contractible.

The normalized Fenchel-Nielsen coordinates map $T_{ls}(X)$ homeomorphically onto $l^{\infty}$ and they map $\overline{T_0(X)}$ homeomorphically onto the space $c_0$ of sequences 
that vanish at infinity. Then $T_{ls}(X)/\overline{T_0(X)}$ is homeomorphic to $l^{\infty}/c_0$ and thus contractible.

To obtain the contractibility for $T_{ls}(X)/\overline{T_{qc}(X)}$, we separate the length and the twist parts of the Fenchel-Nielsen coordinates. Note that the length parts for both $T_{ls}(X)$ and $\overline{T_{qc}(X)}$ cover $l^{\infty}$ and that the quotient of normalized twist parts is $l^{\infty}/c_0$ and thus contractible.
{\it End of proof of Corollary 1.}

\section{The lengths of closed geodesics under quasiconformal maps}

Wolpert \cite{Wol1} proved that a $K$-quasiconformal map distorts the lengths of simple closed geodesics by at most multiplying them with $K$ or at least multiplying them with $1/K$. We prove that this estimate is not sharp for infinite surfaces. To do so, we find a sequence of simple closed curves with lengths going to infinity such that a $K$-quasiconformal map changes their lengths by at most an additive constant.

\vskip .2 cm

\noindent {\it Proof of Proposition 1.}  
Let $l_{\alpha_n}(X)>0$ and
let $\gamma_{n}$ be a shortest simple closed geodesic intersecting $\alpha_{n}$ as before. Then $\gamma_{n}$ is contained in the union of two geodesic pairs of pants $P_{n}^1$ and $P_{n}^2$ of the pants decomposition $\mathcal{P}=\{\alpha_n\}$ of $X$. 

Assume first that $P_{n}^1\neq P_{n}^2$.
Divide each pair of pants $P_n^i$ into two congruent right angled hexagons $\Sigma_{n,1}^i$ and $\Sigma_{n,2}^i$ by three simple mutually disjoint geodesic arcs orthogonal to pairs of boundary geodesics of $P_n^i$ (cf. Figure 1). If a boundary component of $P_n^i$ is a puncture then arcs ending at the puncture have infinite length and there is no requirement to be orthogonal at this boundary component. 

\begin{figure}
\centering
\includegraphics[scale=0.5]{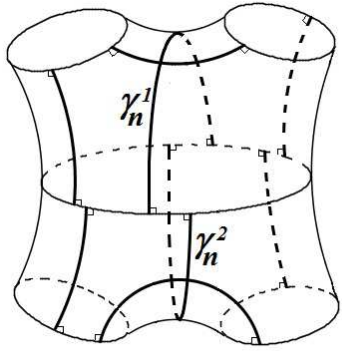}
\caption{Estimating the length of $\gamma_n$.}
\end{figure}

Let $\gamma_{n,1}^i$ be the geodesic arc perpendicular to the side of $\Sigma_{n,1}^i$ contained in $\alpha_n$ and to the opposite side of the hexagon $\Sigma_{n,1}^i$. Define 
$\gamma_{n,2}^i$ to be the geodesic arc perpendicular to the side of $\Sigma_{n,2}^i$ contained in $\alpha_n$ and the opposite boundary side of $\Sigma_{n,2}^i$.
Then $\gamma_{n}^i=\gamma_{n,1}^i\cup \gamma_{n,2}^i$ is a geodesic arc in $P_{n}^i$ with both endpoints at $\alpha_{n}$ that is also orthogonal to $\alpha_{n}$ at these points (cf. Figure 1). 

The geodesic arc $\gamma_{n,1}^i$ divides $\Sigma_{n,1}^i$ into two right angled pentagons. The hexagon $\Sigma_{n,1}^i$ has one boundary side equal to half of the geodesic $\alpha_n$. It follows that at least one of the two pentagons has boundary side $a$ adjacent to $\gamma_{n,1}^i$ equal to at least $1/4$ of the geodesic $\alpha_n$ (and less than $1/2$ of $\alpha_n$). Note that this pentagon has $1/2$ of a geodesic $\alpha_{n'}$ in $\mathcal{P}$ on its boundary and that the length of this boundary side is bounded by $\frac{P}{2}$ from the above. Then \cite{Beardon} we have
$$
\cosh \frac{l_{\alpha_{n'}}(X)}{2}=\sinh l_{\gamma_{n,1}^i}(X)\sinh a.
$$

Since $\frac{\sinh a}{a}\to 1$ as $a\to 0$, it follows that there exists $C_1,C_2>0$ depending on such that
$$
C_1a\leq\sinh a\leq C_2a.
$$
Moreover, we have
$$
1\leq \cosh \frac{l_{\alpha_{n'}}(X)}{2}\leq \cosh \frac{P}{2}=C_3.
$$
The above two inequalities imply that there exist $C_4,C_5>0$ such that
$$
C_4\leq e^{l_{\gamma_{n,1}^i}(X)} l_{\alpha_n}(X)\leq C_5
$$
which implies that for some $C_6>0$
$$
|l_{\gamma_{n,1}^i}(X)-|\log l_{\alpha_n}||\leq C_6.
$$

Thus the length $l_{\gamma_{n}^i}(X)$ is up to a bounded additive constant equal to $2|\log l_{\alpha_{n}}(X)|$ for $i=1,2$.
The geodesic $\gamma_{n}$ is homotopic to the union of $\gamma_{n}^1$ and $\gamma_{n}^2$ and two sub-arcs of $\alpha_{n}$. It follows that
there exist two constants $C_7,C_8\in\mathbb{R}$ such that, for all $n$,
$$
C_7+4|\log l_{\alpha_{n}}(X)|\leq l_{\gamma_{n}}(X)\leq C_8+4|\log l_{\alpha_{n}}(X)|.
$$ 

A $K$-quasiconformal map $f:X\to Y$ changes the lengths of $l_{\alpha_{n}}(X)$ by at most a multiplicative constant between $1/K$ and $K$. However, this only changes the additive constants in the above inequality. Since $l_{\alpha_{n_k}}(X)\to 0$ as $n_k\to 0$, we get that $|l_{\gamma_{n}}(Y)-l_{\gamma_{n}}(X)|\leq C$ for some $C>0$. 

Assume next that $P_{n}^1= P_{n}^2=P_n$. Let $\gamma_{n}^{*}$ be the simple 
geodesic arc in $P_n$ which is orthogonal to $\alpha_n$ at both endpoints.
Then $P_n$ is divided into two hexagons by drawing two more simple mutually  non-intersecting geodesic arcs orthogonal to $\alpha_n$ and the other boundary component $\alpha_{n'}$ of $P_n$. Fix one hexagon $\Sigma_n$ and denote by $\alpha_{n,1}$ and $\alpha_{n,2}$ the two boundary sides of $\Sigma_n$ that lie on $\alpha_n$ (and therefore are adjacent to $\gamma_n^{*}$). Denote by $\alpha_{n,3}$ the boundary side of $\Sigma_n$ on $\alpha_{n'}$ (and therefore opposite to $\gamma_n^{*}$). 
Then by \cite{Beardon} we have
$$
\cosh l_{\gamma_n^{*}}(X)\sinh l_{\alpha_{n,1}}(X)\sinh l_{\alpha_{n,2}}(X) =
\cosh l_{\alpha_{n,3}}(X)+\cosh l_{\alpha_{n,1}}(X)\cosh l_{\alpha_{n,2}}(X) .
$$
Since $l_{\alpha_{n,1}}(X)=l_{\alpha_{n,2}}(X)=\frac{l_{\alpha_n}(X)}{2}\leq P/2$ and $l_{\alpha_{n,3}}(X)\leq P$, the above equality implies that there exist $C_1,C_2>0$ such that
$$
C_1\leq e^{l_{\gamma_n^{*}}(X)}l_{\alpha_n}(X)^2\leq C_2
$$
which implies the result as in the previous case.
{\it End of proof of Proposition 1.}

\end{document}